\documentclass[10pt]{article}
\usepackage{epsfig}
\input epsf
\usepackage{latexsym}
\usepackage{bbold}
\usepackage{graphicx,color}
\usepackage{amsmath, amssymb, graphics}
\usepackage{graphicx}
\newtheorem{theorem}{Theorem}[section]

\newtheorem{cor}{Corollary}[section]

\newtheorem{lemma}{Lemma}[section]

\title{  Anisotropic Wavefronts and Laguerre Geometry  }
\author{  B{\footnotesize ENNETT} P{\footnotesize ALMER}}
\begin{document}
\maketitle
\begin{abstract} 
Motivated by the study of wave fronts in anisotropic media, we propose an incidence geometry of anisotropic spheres in a Finsler-Minkowski space.
An anisotropic version of the Laguerre functional is considered.  In some circumstances, this functional can be used to determine that two wavefronts observed at distinct times in a homogeneous, anisotropic medium, do not originate from the same source.
\end{abstract}

\section{Introduction}

A surface $X:\Sigma \rightarrow {\bf R}^3$ can be regarded as a  source for rays traveling in an isotropic medium with unit velocity. Huygen's Principle tells us that to compute wave front at time $t$ in the future, we can regard each point of the original surface as a point source  of a spherical wave and then take the envelope of the resulting set of spherical waves at time $t$. The resulting wave front is then  given by the parallel surface $X+t\nu$, where $\nu$ is the unit normal to the surface $X$. 

   In general, wave velocity  depends on  the direction in which the wave is traveling and may  depend on position as well. A given material may be isotropic for some types of waves and anisotropic for others, e.g. acoustically isotropic but optically anisotropic.  Here, we will only consider the case where wave velocity is dependent on direction but independent of position.  In this case, the wave fronts for any point source are given by  rescalings of  a fixed shape $W$ which, under reasonable physical assumptions, can be assumed to be  convex \cite{GF}. Because of this, we can  consider $W$ to be the unit sphere of a fixed norm $T$ on ${\bf R}^3$.   Huygen's Principle, which is based on Fermat's Principle,  still holds \cite{A}; if a surface is regarded as a source for waves traveling with constant directionally dependent normal velocity, then the future wave front at time $t$ is an envelope of wave fronts emanating from all point sources in the original surface.  Thus the wave fronts are thus  given by $X+t\xi$, where $\xi$ is the Cahn-Hoffman field of the surface $X$. For  example, for light waves in crystals, the wave fronts for the extraordinary waves with point sources are ellipsoidal For seismic waves propagating in a crystalline solid, the wave front of a point source is a super-ellipse \cite{YN}. 
 
 Laguerre geometry is a sub geometry of the Lie sphere geometry.   
In Laguerre geometry, a surface is embedded in the space of light rays emanating from the surface, each ray being considered as null lines in Lorentz Minkowski space. For each $p\in \Sigma$, this null line is given by $t\rightarrow (X(p),0)+t(\nu(p),1)$.  An element $\alpha$ of the orthogonal group $O(4,1)$ maps the set of null lines emanating from the surface  $X(\Sigma)$ into the set of null lines emanating from a new surface $X'(\Sigma)$, thus $\alpha$ performs a transformation of surfaces as well. The principal aim of Laguerre geometry, which was developed by Blaschke and his followers at the start of the twentieth century \cite{BIII}, is to investigate the invariants of the surface with respect to this action of $O(4,1)$.  In recent times, Laguerre geometry has found applications to ray tracing in computer aided design \cite{PP}.

In this note, we will develop the basics of a formalism which uses four dimensional space to represent the space of anisotropic spheres for a smooth norm in ${\bf R}^3$. Our approach uses a type of Lorentzian Finsler metric, called a conical Finsler metric, which was recently introduced by Javaloyes and Sanchez \cite{JS}.  This allows for the representation of oriented anisotropic spheres, which are the wavefronts having point sources,  as points in  ${\bf R}^4$, in such a way that the incidence relation, the fact that that a point lies on a sphere, is expressible as a homogeneous  equation. The totality  of anisotropic spheres with centers on a surface defines a real line bundle over the surface. A canonical section of this bundle is found which is analogous to the  middle sphere congruence in classical Laguerre geometry. The area of this congruence is used to define an anisotropic Laguerre functional. This functional can be used to detect when two wafefronts, observed at distinct times, do not come from the same source.

By representing a surface using the inverse of its Gauss map, the Euler-Lagrange equation for the Laguerre functional is considered and, as in the isotropic case, the Euler- Lagrange equation  is a linear fourth order equation in an appropriate gauge.

   We would like to thank Professor Dan Dale for  helpful conversations during the preparation of this paper.
\section{Preliminaries}

We will consider wave propagation in an anisotropic, homogeneous material. We do not assume the waves to be of a particular type, it is only assumed that Fermat's Principle holds.  According to this principle,  the path $\alpha $ taken by a ray to reach a point $x$ from a point $y$ minimizes the travel time
among all rays connecting these points, so it  must satisfy the variational principle
$$\delta \int_\alpha \frac{ds}{v(\alpha '(s))}=0\:,$$
where  $v(n)$ denotes the ray speed in the unit direction $n$.

The function $v(\cdot)$ should satisfy reasonable physical assumptions so that $x\rightarrow 1/v(x)$ defines a norm (\cite{GF} Appendix 1), so in the present case, the minimizing path is the line segment from $y$ to $x$ and the minimum travel  time defines a function
$$T_y(x):=\frac{|x-y|}{v(\frac{x-y}{|x-y|})}=\frac{1}{v(x-y)}=T(x-y)\:,$$
where $T=T_0$.
The level set of this function 
$$W_y(t):=\{T_y(x)=t\}\:,$$
is the wave front of rays originating from the point source $y$ at time $t$. To distinguish this from more general wavefronts, we will refer to it as the anisotropic shere with center $y$ and radius $t$. Clearly the anisotropic spheres are all rescalings of a fixed shape $W$ which is the 
unit sphere of the norm $x\rightarrow T(x)$, the so called Wulff shape. 

The gradient $\nabla T_y$ is referred to as the {\it slowness vector}, since the reciprocal of its magnitude $|\nabla T_y|^{-1}$ gives the phase speed i.e. the speed with which the wavefront moves forward with time in the direction normal to the wavefront.  Note that in the anisotropic case, the direction of 
the wavefront at a point does not, in general, coincide with the direction of rays reaching that point.

The dual norm of $T$, denoted $T^*$, is defined for $u\in ({\bf R}^3)^*$ by
$$T^*(u)=\sup_{T(\xi)=1}\langle u,\xi\rangle\:.$$
($\langle ,\rangle$ is the pairing between ${\bf R}^3$ and its dual.)  In the Hamiltonian model of wave propagation, $T^*$ is the Hamiltonian function.
For each $\xi\in W$, there is a unique $\xi^*\in W^*:=\{u\:|\:T^*(u)=1\}$, satisfying $\langle \xi, \xi^*\rangle =1$, $\langle d\xi, \xi^*\rangle =0$.
In addition, there holds \cite{G},
\begin{equation}
\label{dual}
dT_\xi(v)=\langle v,\xi^*\rangle\:,
\end{equation}
for $v\in {\bf R}^3$.

 We associate to the norm $T$  norm a {\it conical Finsler metric}, \cite{JS}, $L$  which is defined for $(x,t)\in {\bf R}^4$ by
$$L((x,t)):=(T(x))^2-t^2\:.$$
Note that $L$ factors as $L=(T-t)(T+t)$ and  $S:=T\mp t$ satisfies the Hamilton-Jacobi equation
$$\partial_{\pm t}S+T^*(\nabla_x S)=0\:,$$
since $\nabla_xS=\nabla T=\xi^*$, by (\ref{dual}). A suitable initial value problem for this equation gives the evolution of a function whose level sets
are an  evolving family of wavefronts whose source is an arbitrary surface \cite{OM}.

\begin{figure}[h!]
   \label{figa}
  \hfill
 \begin{center}
       \framebox{\includegraphics[width=60mm,height=50mm,angle=0]{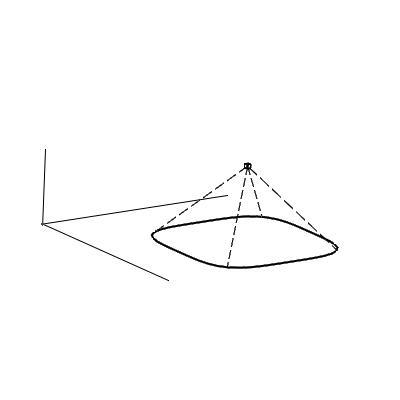}}
    \end{center}
    \caption{An anisotropic sphere represented by a point in ${\bf R}^4$.}
    \end{figure}

\begin{figure}[h!]
   \label{figb}
  \hfill
 \begin{center}
       \framebox{\includegraphics[width=60mm,height=50mm,angle=0]{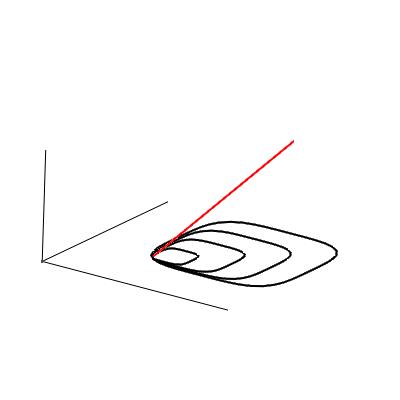}}
    \end{center}
    \caption{The set of anisotropic spheres represented by points lying on a null ray with tangent $(\xi,1)$.}
    \end{figure}

This conical Finsler metric $T$  allows us to regard  ${\bf R}^4$ as the space of oriented anisotropic spheres in ${\bf R}^3$, which, from now on, we identify with  $\{x\in {\bf R}^4\:|\:x_4=0\}$. If we  use $Y=(y,t)$ to represent wavefront $W_t(y)$, (see Figure \ref{figa}), then
we have 
$$W_y(t)=\{X\in {\bf R}^3\:|\: L(Y-X)=0\}\:,$$
where $Y=(y,t)$, $X=(x,0)$. We make the convention that for $t>0$ (resp. $t<0$) the anisotropic sphere is oriented by its outer (resp. inner) normal.
Anisotropic spheres with $t=0$ are called point anisotropic spheres and carry no orientation.
Of course a null line $T_y^2-t^2=0$ in ${\bf R}^4$ will represent an evolving family of anisotropic sphere all passing through the point $y$ in ${\bf R}^3$ which is exactly what one encounters in Huygen's Principle,( see Figure \ref{figb}).  Null vectors will be of the form $(\xi, 1)$, $\xi\in W$ since they satisfy $L((\xi,1))=0$.
    
\section{Anisotropic sphere congruences}

Let $X:\Sigma \rightarrow {\bf R}^3$ be a smooth surface which we regard as either a ray source or as a wavefront. Relative to the norm $T$, we can define an anisotropic normal, known as the {\it Cahn-Hoffman field} \cite{CH}, 
$$\xi:\Sigma \rightarrow W\:$$
defined by the conditions that $\langle dX, \xi^*\rangle \equiv 0$ and that $dX(e_1),dX(e_2),\xi(p)$ is a positively oriented basis of ${\mathbb T}_{X(p)}{\bf R}^3$ whenever $e_1, e_2$ is a positively oriented basis of $\Sigma$.   (We use ${\mathbb T}$ for the tangent space to avoid confusion with the norm $T$.)   
By an {\it anisotropic sphere congruence}, we will mean a smooth map $Y:\Sigma\rightarrow {\bf R}^4$ where $\Sigma$ is a smooth surface. This is just a two parameter family of oriented anisotropic spheres.
An {\it envelope} of $Y$, will be a smooth map $X:\Sigma \rightarrow {\bf R}^3\approx \{x\in {\bf R}^4\:\vert \: x_4=0\}$,  satisfying the conditions:
\begin{eqnarray}
\label{c}
&(i)&\quad L(Y-X )\equiv 0\:,\\
&(ii)&\quad dL_{Y-X} \circ dX\equiv 0\:.\nonumber \end{eqnarray}
These two equations can be interpreted as meaning that  for all  $p\in \Sigma$ the point $X(p)$ lies on the anisotropic sphere represented by $Y(p)$ to first order.
\begin{lemma}
\label{3.1}
For any sufficiently smooth function $\rho$ on $\Sigma$
\begin{equation}
\label{sce} Y:=X+\rho(\xi,1)\:\end{equation}
is an anisotropic spherical congruence enveloped by $X$. Further, (\ref{sce}) is the most general anisotropic spherical congruence enveloped by $X$. 
\end{lemma}
{\it Proof}. For $Y$ defined as above, $Y-X=\rho(\xi,1)$ which is null for $L$ because $T$ is positively homogeneous of degree one.

$$dL=2T dT -2tdt\:,$$
$$dL_{Y-X}=2T(\rho\xi)dT_{\rho \xi}-2\rho dt\:,$$
$$dL_{Y-X}\circ dX=2\rho dT_{\rho \xi}\circ dX\:.$$
Since $T$ is homogeneous of degree one, $dT$ is homogeneous of degree zero, so by (\ref{dual}),  $dT_{\rho \xi}=dT_\xi=\langle \cdot , \xi^*\rangle$. Then
(ii) follows from  the equation $\langle dX, \xi^*\rangle =0$. 

To see that (\ref{sce}) defines the most general anisotropic spherical congruence enveloped by $X$, first note that by (i), it must be that $Y-X$ is null, so we may write
$Y-X=\rho \eta$ with $\eta\in W$. By the argument given above, the equation (ii) would then yield $\langle dX, \eta^*\rangle \equiv 0$. However, there are exactly two points
in $W$ satisfying this equation and those points are $\pm \xi$ where $\xi$ is the Cahn-Hoffman field of $X$. Since only one choice  of these two spheres is consistent with the orientation, the result follows.
{\bf q.e.d.}\\[2mm]

 Canal surfaces are envelopes of one parameter families of anisotropic spheres. Let $\alpha:I\rightarrow {\bf R}^4_1$ denote the curve of anisotropic spheres which we assume to be space-like, i.e. $L(\alpha')>0$. We may assume $\alpha$ is parameterized by arc length and we write $\alpha=({\vec \alpha}, \alpha_4)$. Using Lemma \ref{3.1}, we write the envelope as
\begin{equation}
\label{X}X=\alpha-\alpha_4(\xi,1)=({\vec \alpha}-\alpha_4\xi,0)\:.\end{equation}
 Note that the use of the factor $\alpha_4$ insures that $X$ lies in ${\bf R}^3$. 
 Equation (\ref{c}) (i) holds automatically, since $\alpha-(X,0)=\alpha_4(\xi, 1)$. 
 Equation (\ref{c}) (ii) becomes
\begin{eqnarray*} 0&=&dL_{\alpha-(X,0)}dX=2TdT_{{\vec \alpha}-X}(d\alpha-d\alpha_4 \xi-\alpha_4d\xi)\\
&=&2T\langle \xi^*, {\vec \alpha}'(s)-\alpha_4'(s)\xi\rangle,
\end{eqnarray*}
so $X$ envelopes $\alpha$ if and only if
\begin{equation}
\label{alpha} \langle \xi^*_{{\vec \alpha}-X}, {\vec \alpha}'(s)\rangle =\alpha_4'(s)\:\end{equation}
holds.

We have
$$T({\vec \alpha}'(s))=\sup_{\xi^*\in W^*}\langle \xi^*, {\vec \alpha'}\rangle >|\alpha_4'(s)|\:$$
since $\alpha$ is space-like. On $W^*$, the function $\xi^*\mapsto \langle {\vec \alpha}'(s), \xi^*\rangle$ assumes every value in the interval $[-T({\vec \alpha}'(s)),T({\vec \alpha}'(s))]$ and every value in the interval $(-T({\vec \alpha}'(s)),T({\vec \alpha}'(s)))$ is regular. It follows that for each $s\in I$, this equation (\ref{alpha}) defines a closed curve  $C(s)\subset  W^*$. Letting $\xi$ range throughout the corresponding curve in $W$ for each $s$ produces the envelope $X$ via (\ref{X}).

\begin{figure}[ht]
\label{fign}
{\includegraphics[width=4cm,height=4cm]{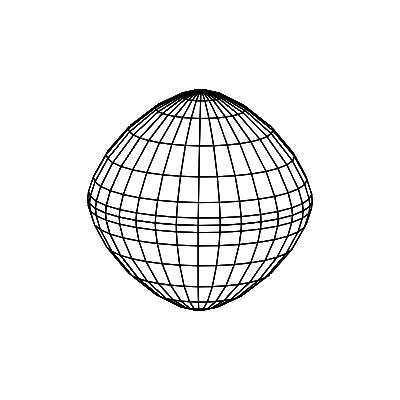}\includegraphics[width=5cm,height=5cm]{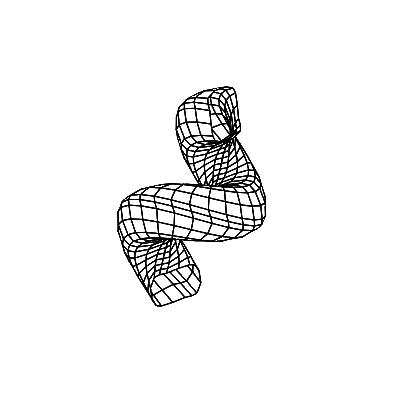}\includegraphics[width=4.5cm,height=4.5cm]{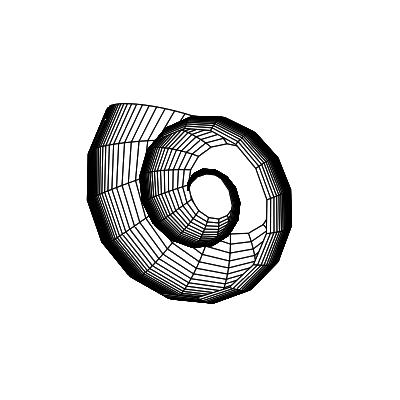}}
\caption{The Wulff shape (left) and two canal surfaces whose source curve is a helix.}
\end{figure}

If $\alpha_4'(s)\equiv 0$, then we can regard the canal surface as a wave front of a disturbance originating on the curve ${\vec \alpha}(s)$ at time $\alpha_4=0$.
A canal surface of this type is shown in the center image of Figure \ref{fign}. Its source curve is a helix in ${\bf R}^3$.  In this case, the surface is an envelope of anisotropic spheres all having the same radius. One can thus regard this surface as a type of tube, although its geometry is more complicated than in the isotropic case since it is not necessarily the product of the curve $\alpha$ with a fixed curve,  The Wulff shape is shown on the left. The final figure shows part of a canal surface which is the envelope of the family of anisotropic spheres $s\mapsto (\cos(s), \sin(s),0,\lambda s)$ so the radii of the anisotropic spheres are not constant.

Huygen's Principle states that given a wavefront $\Sigma$ at time $t_0$, the future wavefront at time $t+t_0$ is an envelope of the set of 
anisotropic spheres of radius $t$ whose centers lie in the surface  $\Sigma$. (See Figure 4.) Huygen's Principle is a direct consequence of Fermat's Priniciple \cite{A}.

With our current formalism, we can prove Huygen's Principle as follows. If a surface in ${\bf R}^3$ is given as the image of an embedding $X:\Sigma\rightarrow {\bf R}^3$, then for each $t>0$, $Y:=(X,t)$ represents the anisotropic sphere with center $X$ and radius $t$.
The parallel surface given by ${\tilde X}:=X+t\xi$ satisfies $Y-{\tilde X}=-t(\xi,-1)$ and so $L(Y-{\tilde X})=0$. Also, by the steps given above and 
(\ref{dual}),
$$dL_{Y-{\tilde X}}\circ d{\tilde X}=\langle \xi^*, d{\tilde X}\rangle=\langle \xi^*, dX+td\xi\rangle \equiv 0\:.$$
We conclude that the for each $t$, the  parallel surface ${\tilde X}$ envelopes $Y$. This means that the secondary wave fronts are given by the maps $X+t\xi$ which is the statement of Huygen's Principle. 
  \begin{figure}[h!]
   \label{C}
  \hfill
 \begin{center}
       \framebox{\includegraphics[width=70mm,height=70mm,angle=0]{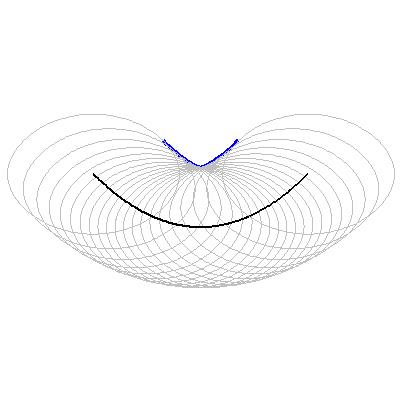}}
    \end{center}
    \caption{The source surface (below) and the wave front (above) as an envelope of anisotropic spheres. In this case, the norm is given by $T(\nu)=1-0.3\nu_3^2$, $\nu \in S^2$.}
    \end{figure}
    
    The manifold ${\bf R}^3\times W$ is a contact manifold with 1-form $\vartheta_{(x,\xi)}:=\langle dx,\xi^*\rangle$. Borrowing terminology from classical differential geometry, we will refer to elements of ${\bf R}^3\times W$ as contact elements. The fiber of the trivial real line bundle 
    $${\bf R}\rightarrow \Lambda \rightarrow {\bf R}^3\times W\:$$
    over a contact element $(x,\xi)$ represents the family of anisotropic spheres passing through $x$ all of which have Cahn-Hoffman field $\xi$ at $x$.
  The tangent space to $\Lambda$ splits,
  $${\mathbb T}_{(x,0)+t(\xi, 1)}\Lambda={\mathbb T}_{(x,\xi)}({\bf R}^3\times W)  \oplus {\bf R}(\xi, 1)\:.$$
  The summands on the right are the horizontal and vertical spaces of the line bundle.

  For a smooth surface $X:\Sigma
    \rightarrow {\bf R}^3$ with Cahn-Hoffman field $\xi$, $f:=(X,\xi):\Sigma \rightarrow {\bf R}^3\times W$ defines a Legendre submanifold , i.e.
    $f^*\vartheta\equiv 0$ holds. Note that an anisotropic spherical congruence is just a section of the bundle $f^*\Lambda$.  For any lift ${\tilde f}$
    of $f$ to $\Lambda$, we let $(d{\tilde f})^T$ denote the projection to the horizontal space of its derivative. The previous lemma states that the lifts 
    are exactly the anisotropic spherical congruences.
    
      Let $X:\Sigma\rightarrow {\bf R}^3$ be an oriented surface with non vanishing curvature. Let $\lambda_i$, $i=1,2$ denote the anisotropic principle curvatures and let $e_i$ denote the corresponding 
   principal directions so that $d\xi(e_i)=-\lambda_idX(e_i)$, $i=1,2$. We normalize the $e_i$'s to have Euclidean norm equal to one,  however, the anisotropic principal directions are not, in general, orthogonal. The {\it anisotropic curvature sphere congruences} are defined by 
   \begin{equation}
   \label{acsc}
   Z_i:=X+\frac{1}{\lambda_i}(\xi, 1)\:, i=1,2\:.\end{equation}
  
   \begin{lemma}
   (cf. \cite{JH}, pg. 61.) Among all anisotropic spherical congruences $Y$ enveloped by $X$, the $Z_i$'s are characterized by the property that for all $p\in \Sigma$, there exists a vector $v\in {\mathbb T}_p\Sigma$ such that 
   $dY(v)$ is null, i.e. $(dY)^T(v) \equiv 0$.\end{lemma}
   {\it Proof.} The proof is immediate; if $Y$ is given by (\ref{sce}), then $(dY)^T(v)=(dX(v)+\rho d\xi(v),0)$ and the result follows. {\bf q.e.d.}\\[4mm]
 We define the {\it anisotropic middle sphere congruence} by
   $$Z:=X+\frac{1}{2}(\frac{1}{\lambda_1}+\frac{1}{\lambda_2})(\xi,1)\:.$$

    Note that the Cahn-Hoffman field $\xi$ has the property that ${\mathbb T}_{\xi(p) }W=dX({\mathbb T}_p\Sigma)$. Therefore we can consider $d\xi_p$ as an endomorphism of $dX({\mathbb T}_p\Sigma)$ to itself and the same is true for $(dY)^T_p$ for any spherical congruence $Y$.
     
         The next result will characterize  $Z$  among all sections of the bundle $\Lambda$. In classical Laguerre geometry, the middle sphere congruence can be characterized as being the unique sphere congruence which defines a conformal map into ${\bf R}^4_1$ if the surface is endowed with its third fundamental form, which is the metric induced by the Gauss map. The analogous property in the anisotropic case would be that the anisotropic middle sphere congruence is conformal with respect to the metric induced by the Cahn-Hoffman map. This will only hold for the class of surfaces for which $d\xi$ is self adjoint, so we look for another characterization.
   \begin{lemma} Among all anisotropic spherical congruences enveloped by $X$, the middle sphere congruence  $Z$ is characterized by the property that  
   $${\rm trace} (d\xi)^{-1}(dZ)^T=0\:.$$\end{lemma}
   {\it Proof.}  For $Y$ given by (\ref{sce}), we compute
   $(dY)^T(e_i)=dX(e_i)+\rho d\xi(e_i)=(\rho-\lambda_i^{-1})d\xi(e_i)\,$
   so the condition that the trace vanishes reduces to the equation $2\rho-(\lambda_1^{-1}+\lambda_2^{-1})=0$. {\bf q.e.d.}\\[2mm]
   \section{The anisotropic Laguerre functional}
   We define the anisotropic Laguerre functional to be the area of $\Sigma$ with respect to the metric $(dZ)^T\cdot (dZ)^T$, that is 
   \begin{equation}
   \label{Lag}{\mathcal L}[X]:= \frac{1}{4}\int_\Sigma (\frac{1}{\lambda_1}-\frac{1}{\lambda_2})^2\:dW\:,\end{equation}
   where $dW$ is the area form of the pull-back metric $d\xi \cdot d\xi$ from $W$. In the isotropic case, this is just the area of the spherical congruence $Z$ in Lorentz-Minkowski space.
   
   Because it was defined as the area of a canonical section of $\Lambda$, it is clear that the functional ${\mathcal L}$ has the same value for all parallel surfaces $X+t\xi$, $t\in {\bf R}$ and should thus be
   considered as a functional of the configuration of null rays originating from $X$. One  application is the following. Suppose that two wavefronts $\Sigma_i$, $i=1,2$ in an anisotropic, homogeneous medium are observed, possibly at distinct times $t_i$, $i=1,2$. Then, one could determine that these two wavefronts did not originate from the same source by showing ${\mathcal L}[\Sigma_1]\ne {\mathcal L}[\Sigma_2]$.  Note that for $|t|>>0$,
{\it  all} wavefronts $X+t\xi$ will start to resemble a rescaling by $|t|$ of $\xi(\Sigma)\subset W$ and that ${\mathcal L}(A)=0$ for all $A\subset W$.

       Since our surfaces are assumed to have non vanishing curvature, we can locally  parameterize them by the inverses of their Gauss maps. This parameterization is global if the surface is convex. It the advantage of expressing
       the entire immersion in terms of one function which is the support function $q$ regarded as a function on the two sphere $S^2$. Then  $X=Dq+q\nu$ where $Df$ denotes the gradient of a function
       $f$ on $S^2$. The derivative of $X$ is expressed $dX=D^2q+qI$ where $D^2$ denotes the Hessian and $I$ is the identity on $TS^2$,(see \cite{E}). If $\tau:=T \vert_{S^2}$, then the analogous parameterization of the Wulff shape $W$ is $\xi=D\tau +\tau \nu$ and $d\xi=D^2\tau +\tau I$. It then follows that $\eta=-(\lambda_1^{-1}+\lambda_2^{-1})$ is given by (\ref{eta}).
   
   In order to state the following result, let $\tau$ be the restriction of the norm $T$ to $S^2$, let $K_W$ denote the curvature
   of the Wulff shape $W$ considered as a function on $S^2$ and let ${\hat \Delta}$ be the Laplacian on $S^2$. If $d\sigma$ is the area form on $S^2$, note that $dW=(1/K_W)\:d\sigma$ holds. 
   \begin{theorem}
 A surface is a critical point of the functional (\ref{Lag}), if and only if the support function $q$ of $\Sigma$, considered as a function on on the 2-sphere, satisfies the fourth order  elliptic linear equation
   \begin{equation}
   \label{pde}  F[q]:=\frac{1}{K_W}{\rm tr}\bigl((D^2\eta+\eta I)(D^2\tau +\tau I)^{-1}\bigr)-2({\hat \Delta} +2)q=0\:,
   \end{equation}
   where 
   \begin{equation}
   \label{eta} \eta:={\rm tr}\bigl((D^2q+qI)(D^2\tau +\tau I)^{-1}\bigr)=-\frac{1}{\lambda_1}-\frac{1}{\lambda_2}\:.\end{equation}

      \end{theorem} 
      {\bf Remark} In the classical (isotropic) case $W=S^2$, $\tau \equiv 1$, $K_W\equiv 1$  and (\ref{pde}) becomes $0=({\hat \Delta}+2)\eta-2({\hat \Delta}+2)q=({\hat \Delta}+2)\eta-2\eta={\hat \Delta}\eta$. This means that
      $\eta =-(k_1^{-1}+k_2^{-1})$ is harmonic on $S^2$ or equivalently that the support function $q$ satisfies $({\hat \Delta}+2){\hat \Delta}q=0$. This equation can be found in the work of Blaschke \cite{B}  and it can be reduced to the biharmonic equation in ${\bf R}^2$ \cite{PG}, using a twistor correspondence.
      
      Also, in the case where $\Sigma =W$, we have $\lambda_1\equiv -1\equiv \lambda_2$ so the integrand in ${\mathcal L}$ vanishes identically, i.e., $W$ (and its rescalings) are strict minima of the functional for any boundary conditions. In this case $q\equiv \tau$ and, as expected,  the equation (\ref{pde}) holds.
      \\[4mm]
  {\it Proof.}  Write  $(\frac{1}{\lambda_1}-\frac{1}{\lambda_2})^2=\eta^2-4/\lambda_1\lambda_2=\eta^2-4$det$((D^2q+qI)(D^2\tau +\tau I)^{-1})$. When we make a compactly supported variation of  the surface through parallel immersions, we have 
     $q\rightarrow q+\epsilon {\dot q}+...$, ${\dot q}\in C^\infty_c$ and $\tau$ remains constant. Therefore $\delta \eta=$tr$((D^2{\dot q}+{\dot q}I)(D^2\tau +\tau I)^{-1}$.  If $A(\epsilon)$ is a smooth one parameter curve of matrices, then $\partial_\epsilon$det$A(\epsilon)_{\epsilon=0}=$det$A(0)$tr$(A^{-1}(0){\dot A}(0))\:$. Using this with $A:=(D^2q+qI)(D^2\tau +\tau I)^{-1}$, we obtain
\begin{eqnarray*} && \delta {\rm det}((D^2q+qI)(D^2\tau +\tau I)^{-1})\\&=& {\rm det}\bigl[(D^2q+qI)(D^2\tau +\tau I)^{-1}\bigr]{\rm tr}\bigl((D^2\tau +\tau I)(D^2q +q I)^{-1}(D^2{\dot q}+{\dot q}I)(D^2\tau +\tau I)^{-1}\bigr)\\
&=&\frac{1}{\lambda_1\lambda_2}{\rm tr}\bigl( (D^2{\dot q}+{\dot q}I)(D^2q+qI)^{-1}\bigr)\:.\end{eqnarray*}
For a variation through parallel immersions, the area form $dW$ is unchanged since the Cahn-Hoffman map factors through the Gauss map $\nu$.  We then have
\begin{eqnarray*}\delta \label{L}{\mathcal L}[X]&=&\int_\Sigma 2\eta \delta \eta \:dW-4\int_\Sigma \delta  {\rm det}((D^2q+qI)(D^2\tau +\tau I)^{-1})\:dW\\
&=&\int_W 2\eta {\rm tr}((D^2{\dot q}+{\dot q}I)(D^2\tau +\tau I)^{-1}\:dW-4\int_\Sigma \frac{1}{\lambda_1\lambda_2}{\rm tr}\bigl( (D^2{\dot q}+{\dot q}I)(D^2q+qI)^{-1}\bigr)\:dW
\end{eqnarray*}
For the first integral, since det $(D^2\tau +\tau I)=1/K_W$ and $dW=d\sigma/K_W$, we get
\begin{eqnarray*}
{\rm tr}((D^2{\dot q}+{\dot q}I)(D^2\tau +\tau I)^{-1}\:dW&=&K_W\bigl(M({\dot q},\tau)+{\dot q}{\hat \Delta}\tau +\tau{\hat \Delta} {\dot q}+2{\dot q}\tau \bigr)\frac{d\sigma}{K_W}\:.\\
&=&\bigl(M({\dot q},\tau)+{\dot q}{\hat \Delta}\tau +\tau{\hat \Delta} {\dot q}+2{\dot q}\tau \bigr)d\sigma\:.
\end{eqnarray*}
With the help of formula (\ref{ibp}) below with $S=S^2$, one can obtain
\begin{eqnarray}
\label{A}
2\int_\Sigma \eta \delta \eta\:dW&=&2\int_\Sigma {\dot q}\bigl(M(\tau, \eta)+\eta({\hat \Delta} \tau +2\tau)+\tau {\hat \Delta}\eta\bigl)\:d\sigma \nonumber \\ 
&=& 2\int_\Sigma \frac{{\dot q}}{K_W}\: {\rm tr}\bigl((D^2\eta+\eta I)(D^2\tau +\tau I)^{-1}\bigr)\:d\sigma\:.
\end{eqnarray}

In the second integral, we have,
\begin{eqnarray}
 \frac{1}{\lambda_1\lambda_2}{\rm tr}\bigl( (D^2{\dot q}+{\dot q}I)(D^2q+qI)^{-1}\bigr)\:dW&=&\frac{K_w}{K_\Sigma}K_\Sigma\bigl(M({\dot q}, q)+{\dot q}({\hat \Delta} q+2q)+q{\hat \Delta} {\dot q}\bigr)\:\frac{d\sigma}{K_w} \nonumber \\ 
 &=&\bigl(M({\dot q}, q)+{\dot q}({\hat \Delta} q+2q)+q{\hat \Delta} {\dot q}\bigr)\:d\sigma\:.
 \end{eqnarray}
 Using (\ref{ibp}) with $w\equiv 1$, one easily obtains

\begin{equation}
\label{B}
-4\int_\Sigma \frac{1}{\lambda_1\lambda_2}{\rm tr}\bigl( (D^2{\dot q}+{\dot q}I)(D^2q+qI)^{-1}\bigr)\:dW=-4\int_\Sigma {\dot q}({\hat \Delta}+2)q\:d\sigma\:.
\end{equation}
 Replacing the corresponding terms in the formula for $\delta{\mathcal L}$ by (\ref{A}) and (\ref{B}) results in the equation (\ref{pde}).

 The principal part of this equation is given by ${\rm tr}\bigl((D^2\eta)(D^2\tau+\tau I)^{-1})$, which is elliptic since the matrix $D^2\tau+\tau I$ is positive definite everywhere on $S^2$ because $\tau$ is the restriction of a norm. {\bf q.e.d.}\\[2mm]
 {\bf Remarks.} As mentioned above, under rescaling of the surface $X\rightarrow rX$,  the functional rescales according to  ${\mathcal L}\rightarrow r^2{\mathcal L} $.
 It follows that ${\mathcal L} $ has no closed critical points other than anisotropic spheres since these are the only closed surfaces for which 
 ${\mathcal L}=0$.( The last statement characterizing anisotropic spheres basically follows from lemma 3.3 of \cite{HH}.) In light of this, boundary conditions should be imposed. Since the Euler-Lagrange equation is fourth order, it seems natural to fix the boundary of the surface to first order.
 
 The calculations in the proof show that the first variation of the function ${\mathcal L}$ is given by the
 $$\delta {\mathcal L}=\int_\Sigma {\dot q}F[q]\:d\sigma\:.$$
 For a critical point, we have $F[q]=0$ and, taking into account the linearity of $F$, we obtain the second variation immediately as,
 $$\delta^2 {\mathcal L}=\int_\Sigma {\dot q}F[{\dot q}]\:d\sigma\:.$$
 Higher order variations can be computed inductively.

 An immediate consequence of the theorem is that there is an abundance of examples of critical surfaces for ${\mathcal L}$. 
 As far as explicit examples go, we have, besides the Wulff shape,  the following. 

\begin{cor} For any norm $T$ having axially symmetric Wulff shape, the helicoids $(r,\theta)\mapsto(r\cos\theta, r\sin \theta, a\theta)$, $a\in {\bf R}^*$,  are critical points of the anisotropic Laguerre functional ${\mathcal L}$.
\end{cor}
 {\it Proof.} First of all, the helicoids have non vanishing curvature. Also, they have the remarkable property that for {\it any} axially symmetric norm, their anisotropic mean curvature $\Lambda:=\lambda_1+\lambda_2$ vanishes. \cite{KP}. Therefore $\eta=\lambda_1^{-1}+\lambda_2^{-1}\equiv 0$ and, in particular,  ${\hat \Delta}q+2q=-(k_1^{-1}+k_2^{-1})\equiv 0$, so (\ref{pde}) holds, i.e. they are minimal in the usual sense. {\bf q.e.d.}
  \section{Other Gauge Invariants}
  
  Here we show how to obtain other gauge invariants of the bundle $f^*\Lambda$, that is, quantities which are the same for all parallel immersions in a family
  $\{X+t\xi \:|\:t\in {\bf R}\}$. At first, we restrict our attention in this section to closed convex surfaces. 
  
Associated with a given Wulff shape $W:=\{x\:|\: T(x)=1\}$ is an anisotropic surface energy functional which assigns to an oriented, immersed surface in ${\bf R}^3$  the valus
$${\mathcal F}[\Sigma]:=\int_\Sigma T^*(\nu)\:d\Sigma\:.$$

  Let ${\mathcal F}(t)$ denote the area of the surface $X+t\xi$. By using the definition of the anisotropic principle directions $dX(e_i)=-1/\lambda_i d\xi(e_i)$,
 the  measure $d{\mathcal F}_p(t)$ for the energy of the surface  $X+t\xi$ at $p\in \Sigma$ is $d{\mathcal F}_p(t):=\pm T^*(\nu)(t-\lambda_1)^{-1}(t-\lambda_2)^{-1}dW_p$, i.e.
  $$d{\mathcal F}_p(t)=\pm T^*(\nu)(t^2+\eta t+ K_W/K_\Sigma)dW_p =:Q(p,t)dW_p\:.$$
Replacing $X$ by $X+t_0\xi$ translates the $1/\lambda_i$'s by $t_0$ so the discriminant $({\rm disc}(Q(p,t))^2$ is unchanged. Note that his quantity is exactly the density $(1/\lambda_1-1/\lambda_2)^2$ appearing in the Laguerre functional, which gives another way to see that it is gauge invariant.

Now observe that we could just as well have integrated and then taken the discriminant to obtain another gauge invariant
\begin{eqnarray*}
I&:=&{\rm disc}(\int_\Sigma T^*(\nu) Q(p,t)dW_p)\\
&=& (\int_\Sigma T^*(\nu)\eta\:dW)^2-4(\int_\Sigma T^*(\nu)dW)(\int_\Sigma T^*(\nu)\frac{K_W}{K_\Sigma}\:dW)\\
&=& (\int_\Sigma T^*(\nu) \eta\:dW)^2-12{\rm vol}(W){\mathcal F}(\Sigma),\\
\end{eqnarray*}
where, in the last step, we have used that the energy of the Wulff shape is equal to three times its enclosed volume.

The quantity $I$ can be computed in a more direct way starting with the expansion of the anisotropic energy density of $X+t\xi$ which gives
$T^*(\nu)\:d\Sigma(t)=T^*(\nu)(t^2 K_\Sigma/K_w-\Lambda t +1)\:d\Sigma$, where $\Lambda=\lambda_1+\lambda_2$ is the anisotropic mean curvature. Integrating this over $\Sigma$ and taking the discriminant gives:
$$I=(\int_\Sigma T^*(\nu) \Lambda\:d\Sigma)^2-4 {\mathcal F}[\Sigma]\int_\Sigma T^*(\nu)\frac{K_\Sigma}{K_W}\:d\Sigma\:,$$
which is valid for {\it any} sufficiently smooth compact surface with or without boundary.

Since $I$ also rescales quadratically, the quantity $I/{\mathcal L}$ is both gauge invariant and scale invariant. This ratio is defined for all closed, convex surfaces other than anisotropic spheres.

Another gauge invariant can be obtained by using the discriminant cubic which gives the expansion for the volume $V(t)$ of the section $X+t\xi$. The integral
of this polynomial is the so called Steiner polynomial which has been widely studied in Convex Analysis.

   \section{Appendix}
   
   We prove an integration by parts formula for a linearized two dimensional Monge Ampere equation which was used above.
   
   On a smooth surface $S$, we let $M[u]$, $u\in C^{\infty}(\Sigma)$, denote the Monge-Ampere operator given by
   $$-2M[u]=|\nabla ^2u|^2-(\Delta u)^2\:.$$
   We also define a symmetric operator by
   $$M(u,w)=\partial_tM[u+tw]_{t=0}=\partial_tM[w+tu]_{t=0}\:.$$
   We wish to show that for $u,w\in C^{\infty}$, $\zeta\in C^{\infty}_c$, there holds
   \begin{equation}
   \label{ibp} \int_S w(M(u,\zeta)-K\nabla u\cdot \nabla \zeta)\:dS=\int_S \zeta (M(u,w)-K\nabla u\cdot \nabla w)\:dS\:.\end{equation}  
   
   To prove this, we start with the Lichnerowicz formula
   $$\nabla\cdot \bigl[\nabla \frac{1}{2}|\nabla u|^2-\Delta u\nabla u\bigr]=-2M[u]+K|\nabla u|^2\:.$$
   By replacing $u\rightarrow u+tw$ in this formula and differentiating with respect to $t$ at $t=0$, we arrive at
   $$\nabla \cdot ([\nabla(\nabla u\cdot \nabla w)-\Delta u\nabla w-\Delta w\nabla u]=-2M[u,w]+2K\nabla u\cdot \nabla w\:.$$
   By repeated use of Green's formulas, we then have
   \begin{eqnarray*}
  &&\int_\Sigma \zeta \nabla\cdot \bigl[ \nabla(\nabla u\cdot \nabla w)-\Delta u \nabla w-\Delta w\nabla u\bigr]\:d\Sigma\\
   &=& -\int_\Sigma \nabla \zeta\cdot \bigr[ \nabla(\nabla u\cdot \nabla w)-\Delta u \nabla w-\Delta w\nabla u\bigr]\:d\Sigma)\\
   &=&  \int_\Sigma  (\nabla u\cdot \nabla w)\Delta \zeta-w\nabla\cdot(\Delta u\nabla \zeta)+\Delta w(\nabla u\cdot \nabla \zeta)\:d\Sigma\\
   &=& \int_\Sigma- w\nabla\cdot (\Delta \zeta\nabla u)-w\nabla \cdot ( \Delta u\nabla \zeta)+w\Delta(\nabla u\cdot \nabla \zeta)\:d\Sigma\\
   &=& \int_\Sigma w \nabla\cdot \bigl[ \nabla(\nabla u\cdot \nabla \zeta)-\Delta u \nabla \zeta-\Delta \zeta\nabla u\bigr]\:d\Sigma\:,
   \end{eqnarray*}
 which proves (\ref{ibp}).


\bigskip
\begin{flushleft}
Bennett P{\footnotesize ALMER}\\
Department of Mathematics\\
Idaho State University\\
Pocatello, ID 83209\\
U.S.A.\\
E-mail: palmbenn@isu.edu
\end{flushleft}
\end{document}